\documentclass[a4paper,12pt,leqno]{article}

\usepackage{geometry} 

\usepackage{tikz-cd}\usepackage{bussproofs}
\usepackage{amssymb,amsmath}
\usepackage{proof}
\usepackage{psvectorian}
\usepackage{comment}
\usepackage{CJKutf8}
\usepackage[utf8]{inputenc}
\usepackage{graphicx}
\usepackage{graphicx}
\graphicspath{ {images/} }
\usepackage{qtree}
\usepackage{mathtools}
\usepackage{hyperref}
\usepackage{listings}
\usepackage{indentfirst}
\usepackage{titlesec}

\DeclareUnicodeCharacter{00B9}{$^{1}$}
\DeclareUnicodeCharacter{00AC}{$\neg$}

\DeclareUnicodeCharacter{03B5}{$\epsilon$}

\DeclareUnicodeCharacter{2200}{$\forall$}

\DeclareUnicodeCharacter{2203}{$\exists$}

\DeclareUnicodeCharacter{2229}{$\cap$}

\DeclareUnicodeCharacter{222A}{$\cup$}

\DeclareUnicodeCharacter{2192}{$\rightarrow$}
\DeclareUnicodeCharacter{2282}{$\subset$}

\DeclareUnicodeCharacter{03C9}{$\omega$}

\DeclareUnicodeCharacter{2218}{$\circ$}

\DeclareUnicodeCharacter{207B}{$^{-1}$}

\numberwithin{equation}{section}

\newtheorem{defn}[equation]{Definition}

\newtheorem{rem}[equation]{Remark}
\newtheorem{exm}[equation]{Example}

\newtheorem{notat}[equation]{Notation}
\newtheorem{newpar}[equation]{}

\newtheorem{xdefn}{Definition.}
\newtheorem{xproposition}{Proposition.}
\newtheorem{xcorollary}{Corollary.}
\newtheorem{xrem}{Remark.}
\newtheorem{xexm}{Example.}
\newtheorem{xlemma}{Lemma.}
\newtheorem{xtheorem}{Theorem.}
\newtheorem{xnotat}{Notation.}
\newtheorem{xnewpar}{\it}
\newtheorem{xproof}{{\it Proof. }}
\newtheorem{xproofof}{{\it Proof}}

\newenvironment{remark}{\begin{rem}\em}{\end{rem}}

\newenvironment{newparagraph*}[1]{\begin{xnewpar}\hspace*{-1.5mm}{#1}. \rm}{\end{xnewpar}}

\newenvironment{definition*}{\begin{xdefn}\em}{\end{xdefn}}
\newenvironment{remark*}{\begin{xrem}\em}{\end{xrem}}
\newenvironment{example*}{\begin{xexm}\em}{\end{xexm}}
\newenvironment{notation*}{\begin{xnotat}\em}{\end{xnotat}}
\newenvironment{proposition*}{\begin{xproposition}}{\end{xproposition}}
\newenvironment{corollary*}{\begin{xcorollary}}{\end{xcorollary}}
\newenvironment{lemma*}{\begin{xlemma}}{\end{xlemma}}
\newenvironment{theorem*}{\begin{xtheorem}}{\end{xtheorem}}

\titleformat*{\section}{\large\bfseries}

\begin{document}

	\title{Introduction to PyLog}
	\author{Clarence Protin}
	
	%	\date{12th of September 2021}
	\maketitle

 \begin{abstract}
 PyLog is a minimal experimental proof assistant based on linearised natural deduction for intuitionistic and classical first-order logic extended with a comprehension operator. PyLog is interesting as a tool to be used in conjunction with other more complex proof assistants and formal mathematics projects (such as Coq and Coq-based projects).
 Proof assistants based on dependent type theory are at once very different and profoundly connected to the one employed by Pylog via the Curry-Howard correspondence. The Tactic system of Coq presents us with a top-down approach to proofs (find  a term inhabiting a given type via backtracking the rules, typability and type-inference being automated) whilst the classical approach of Pylog follows how mathematical proofs are usually written.
 
 Pylog should be further developed along the lines of Coq in particular through the introduction of many "micro-automatisations"  and a nice IDE.
 \end{abstract}

 \section*{Introduction}
 
 As Voevodsky pointed out in a \href{https://www.youtube.com/watch?v=E9RiR9AcXeE}{public lecture} in Princeton,
 it is highly desirable to obtain a foundations of mathematics that will allow automatic verification of proofs. Standard approaches employ intuitionistic higher order logic or various powerful dependent type theories  that have the added bonus of
 the constructive computational information furnished by the  Curry-Howard isomorphism. Notable examples
 are formal mathematics projects based on Coq and Agda and those based on Homotopy Type Theory.
 PyLog is a computationally and philosophically alternative approach which aims at fulfilling a number of desiderata:
 \begin{enumerate}
 	\item The user environment must be easy, simple, intuitive and attractive to use for the  logician or mathematician and be transparent to the programmer so as to easily facilitate access to data structures and algorithms for future development and applications.
 	\item The process of writing proofs (and the checking algorithm) should be agreable and resemble structurally actual mathematical practice - or at least the template laid down by the Principia Mathematica (1910).
 	\item It should be first-order (with a weak "parametric" second-order extension) with all its type-free simplicity and versatility.
 	\item It should be easy to combine different formalised theories and to organise theorems into theories.
 	\item Formalised proofs  can be easily checked either by a human or a machine.
 	\item The difficulty of formalising and checking a given theory should not exceed mathematical difficulty of the theory involved.
 	\item Classical logic is to be seen as an extension of intuitionistic logic and the user is free to use the classical negation rule or not.
 	
 \end{enumerate}

 The key ingredients that I propose are:
 
 \begin{itemize}
 	\item A linearised natural deduction for the predicate calculus with equality (and some weak second-order extension) extended with
 	a Kelley-Morse style extension operator.
 	\item Can easily formalise Kelley-Morse set theory but is not restricted to it.
 \end{itemize}

 \section*{The Logic of PyLog}
 
 PyLog is based on the natural deduction presentation of first-order predicate logic with equality endowed with rules for a class-forming operator $\{x: P(x)\}$ and a primitive binary predicate $\in$.
 Pylog also includes second-order variables  allowing  us to instantiate logical validities.
 The language of PyLog consists of finite sets of constants, (first-order) variables, second-order variables,  function symbols of different arities $n >0$, predicate symbols of different arities $n<0$ and the special symbol $\bot$.
 Terms and formulas are defined by mutual recursion:

 \begin{itemize}
 \item A constant $c$ is a term.
 \item A variable $x$ is a term.
 \item If $t_1,...t_n$ are terms and $f$ is a $n$-ary function symbol then $f(t_1,...,t_n)$ is a term.
 \item If $A$ is a formula and $x$ is a variable then $\{x : A\}$ is a term (called an \emph{extension})
 \item If $P$ is a $n$-ary predicate symbol and $t_1,...t_n$ are terms then $P(t_1,...,t_n)$is a formula.
 \item If $t$ and $s$ are terms then $t=s$ is a formula (this is a particular case of the last condition).
 \item If $\mathfrak{A}$ is a second-order variable then it is a formula.
 \item If $A$ and $B$ are formulas then $A\vee B$, $A\enspace \&\enspace B$, $A\rightarrow B$ are formulas.
 \item If $x$ is a variable and $A$ is a formula then $\forall x.A$ and $\exists x.A$ are formulas.
 \item $\bot$ is a formula.
 \end{itemize}

We define the set $FV(e)$ of \emph{free variables} of an expression $e$ (term or formula) as follows:

\begin{itemize}
	\item $FV(c) = \emptyset$
	\item $FV(x) = \{x\}$
	\item $FV(f(t_1,...,t_n)) = \bigcup_{i=1,...,n} FV(t_i)$
		\item $FV(P(t_1,...,t_n)) = \bigcup_{i=1,...,n} FV(t_i)$
	\item $FV(\bot) = \emptyset$
	\item $FV(\mathfrak{A}) = \emptyset$
	\item $FV(A\vee B)$, $FV(A\enspace\&\enspace B)$ and  $FV(A\rightarrow B)$ are equal to $FV(A) \cup FV(B)$
   \item $FV(\forall x.A)$, $FV(\exists x.A)$ and $FV(\{x: A\})$ are equal to $FV(A)/\{x\}$ 
\end{itemize}

As usual we consider expression \emph{modulo} the renaming of quantified variables or variables within the scope of an extension: in any subexpression of the form $\forall x.A$, $\exists x.A$ or $\{x:A\}$ we may rename $x$ and all free occurrences of $x$ in $A$ to a fresh variable $y$ as long $y$ does not occur within the scope of some quantifier  $\forall y$, $\exists y$ or extension $\{y:...\}$.
When we write $A[t/x]$ we assume that the bound variables of $A$ have been renamed so as to be distinct from $FV(t)$ (this is a slightly stronger condition than we actually need).

In PyLog proofs are always in the context of a \emph{proof environment}. This consisting of:
 
 \begin{itemize}
 	
  \item A list formulas called \emph{axioms}
  \item A list of formulas called \emph{assumed theorems}
  \item A list  of \emph{defining equations} for constants or function symbols of the form $c = A$ or $f(x_1,...,x_n) =A$
  \item A list of \emph{predicate definitions} consisting of triples $(P, (x_1,...,x_n),A)$ defining $n$-ary predicate symbols $P$. Here
  $FV(A) \subseteq \{x_1,...,x_n\}$. Triples are also denoted by $P(x_1,...,x_n)\equiv A$.

  \end{itemize}

All lists above may be empty. In PyLog the default proof environment consists of a single predicate definition  for $Set(x)$ defined as $\exists y. x \in y$. The language is endowed further with the  primitive binary predicate $\in$. There are no other functions, predicates or constants.

The proof system of PyLog is based on a linearised variant of natural deduction with conservative second-order order extension.

 We first present the system in the standard form. We assume the reader is familiar with proof trees and the concept of dependency(the first chapters of \cite{dag} are sufficient).  The proof system of PyLog consists of the following.
 We have \emph{purely logical rules} which are the rules for \emph{minimal predicate calculus} plus the
 intuitionistic and classical negation rules:

 \begin{center}
 
 \AxiomC{$A$ } 
 \AxiomC{$B$}
 \RightLabel{AndInt}
 \BinaryInfC{$A\enspace\&\enspace B$}
 \DisplayProof
 \hskip 1.5em
 \AxiomC{$A\enspace\&\enspace B$} 
 \RightLabel{AndElimL}
 \UnaryInfC{$A$}
 \DisplayProof
 \hskip 1.5em
  \AxiomC{$A\enspace\&\enspace B$} 
  \RightLabel{AndElimR}
  \UnaryInfC{$B$}
 \DisplayProof
 \end{center}
 
 \begin{center}
 \AxiomC{$[A]$}
 \noLine
 \UnaryInfC{$B$} 
 \RightLabel{ImpInt}
 \UnaryInfC{$A\rightarrow B$}
 \DisplayProof
 \hskip 1.5em
  \AxiomC{$A$ } 
 \AxiomC{$A\rightarrow B$}
 \RightLabel{ImpElim}
 \BinaryInfC{$B$}
 \DisplayProof
 \vspace{5mm}
 
 \end{center}

\begin{center}
 \AxiomC{$A$} 
 \RightLabel{OrIntL}
 \UnaryInfC{$B \vee A$}
  \DisplayProof
 \hskip 1.5em
 \AxiomC{$A$} 
 \RightLabel{OrIntR}
 \UnaryInfC{$A \vee B$}
  \DisplayProof
  \hskip 1.5em
  \AxiomC{$A \vee B$}
  \AxiomC{$[A]$}
  \noLine
  \UnaryInfC{$C$}
  \AxiomC{$[B]$}
  \noLine
  \UnaryInfC{$C$}
  \RightLabel{OrElim}
  \TrinaryInfC{$C$}
   \DisplayProof
 \end{center}
  
 \begin{center}
\AxiomC{$A$} 
\RightLabel{ForallInt$_y$}
\UnaryInfC{$\forall x. A[x/y]$}
\DisplayProof
 \hskip 1.5em
  \AxiomC{$\forall x. A$} 
 \RightLabel{ForallElim}
 \UnaryInfC{$A[t/x]$}
 \DisplayProof
 \end{center}

 \begin{center}
 \AxiomC{$A[t/x]$} 
 \RightLabel{ExistsInt}
 \UnaryInfC{$\exists x. A$}
 \DisplayProof
 \hskip 1.5em
 \AxiomC{$\exists x. A$}
 \AxiomC{$[A[y/x]]$}
 \noLine
 \UnaryInfC{$C$}
 \RightLabel{ExistsElim}
 \BinaryInfC{$C$}
 \DisplayProof
 \end{center}
 
 \begin{center}
 \AxiomC{$\bot$} 
 \RightLabel{Abs$_i$}
 \UnaryInfC{$A$}
 \DisplayProof
 \hskip 1.5em
 \AxiomC{$[\sim A]$}
 \noLine
 \UnaryInfC{$\bot$} 
 \RightLabel{Abs$_c$}
 \UnaryInfC{$A$}
 \DisplayProof
 \end{center}

The proviso for ForallInt is that $y$ cannot occur in any assumption on which $A$ depends and the proviso for ExistsElim is that $y$ does not occur in $\exists y.A$ or in $C$ or on any hypothesis on which $C$ depends other than $[A[y/x]]$. $\sim A$ is syntactic sugar for $A\rightarrow \bot$.

We have the \emph{class rules}:
 \begin{center}
 \AxiomC{$ t \in \{x: A\}$ } 
 \RightLabel{ClassElim}
 \UnaryInfC{$Set(t) \enspace\&\enspace A[t/x]$}
 \DisplayProof
 \hskip 1.5em
 \AxiomC{$Set(t) \enspace\&\enspace A[t/x]$ } 
 \RightLabel{ClassInt$_x$}
 \UnaryInfC{$ t \in \{x: A(x)\}$}
 \DisplayProof
 \end{center}
 
 These rules express the \emph{classification axiom scheme} of Kelley-Morse set theory such as formulated in the appendix of \cite{kel}.
 
 We have also the \emph{equality rules}\footnote{this is inspired by the treatment in \cite{tro}}:

 \begin{center}
 	\AxiomC{}
 	\RightLabel{Identity}
 	\UnaryInfC{$t=t$}
 	\DisplayProof
 	\hskip 1.5em
 	\AxiomC{$s=t$}
 	\RightLabel{Symmetry}
 	\UnaryInfC{$t=s$}
 	\DisplayProof
 	\hskip 1.5em
 	\AxiomC{$A$ } 
 	\AxiomC{$t = s$}
 	\RightLabel{EqualitySub}
 	\BinaryInfC{$A'$}
 	\DisplayProof

 \end{center}
 
 where in EqualitySub $A'$ is $A$ when a specified number of occurrences of $t$ are replace by $s$.
 EqualitySub is of fundamental importance in using defined constants and function symbols of the proof environment.
 
 We have then our \emph{second-order rule}:
 \begin{center}
 	\AxiomC{$A$ } 
 	\RightLabel{PolySub$_{\mathfrak{A}}$}
 	\UnaryInfC{$A'$}
 	\DisplayProof
 \end{center}
 
 where $A'$ results from $A$ by substituting all occurrences of $\mathfrak{A}$ in $A$ by a formula $B$. The proviso is that $\mathfrak{A}$  does not occur in any hypothesis on which $A$ depends and that no free variable in $B$ becomes bound after the substitution.
 
 \begin{remark}
 This rule is to be understood as a combination of an invisible second-order generalisation of the variable $\mathfrak{A}$ followed by an instantiation by $B$.	
 	\end{remark}
 
 The final set of rules concern how information in the proof environment is introduced into the proof.

 \begin{center}
 	\AxiomC{}
 	\RightLabel{AxInt$_n$}
 	\UnaryInfC{$A$}
 	\DisplayProof
 	\hskip 1.5em
 		\AxiomC{}
 	\RightLabel{TheoremInt$_n$}
 	\UnaryInfC{$A$}
 	\DisplayProof
 	\hskip 1.5em
 	 	\AxiomC{}
 	\RightLabel{DefEqInt$_n$}
 	\UnaryInfC{$t = s$}
 	\DisplayProof
 	\hskip 1.5em
 	\AxiomC{$A$ } 
 	\RightLabel{DefExp}
 	\UnaryInfC{$A'$}
 	\DisplayProof
 		\hskip 1.5em
 	\AxiomC{$A$ } 
 	\RightLabel{DefSub}
 	\UnaryInfC{$A'$}
 	\DisplayProof

 \end{center}
 
 The first three rules simply add the $n$th formula in the lists of axioms, assumed theorems and defining equations respectively. DefExp does the following. Assume we have a definition $P(x_1,...,x_n)\equiv B$. Then DefExp replaces specified occurrences of subformulas of the form $P(t_1,...,t_n)$ by $B[t_1/x_1,...,t_n/x_n]$ (the resulting expression is denoted by $A'$ in the rule). DefSub does the inverse of this. For chosen $t_1,...,t_n$ we  must specify the occurrences of expressions of the form $B[t_1/x_1,...,t_n/x_n]$ in $A$ which we wish to "collapse" into $P(t_1,...,t_n)$\footnote{For example if we had the definition $Set(x) \leftrightarrow \exists y. x\in y$ and a line in our proof such as

 \[1. Set(x)\]
 
 then calling DefExp(n,"Set",[0]) would add the line
 
 \[2. \exists y. x \in y\]
 
 On the other hand if we had a line
 
 \[3. \exists y. x\in y \enspace \& \enspace \exists z. y \in z\]
 
 then the command DefSub(3,"Set", ["x"],[0]) would yield
 
 \[4. Set(x)  \enspace \& \enspace \exists z. y \in z\]}.

 When we wish to use defined functions or constants we first introduce the defining equalities into our proof by means o DefEqInt$_n$  and then make use of EqSub.
 The above are the core rules of PyLog.  We introduce the usual abbreviation $A\leftrightarrow B$
 and so finally have two rules to toggle this notation:

 \begin{center}
 	\AxiomC{$A\rightarrow B \enspace\&\enspace B \rightarrow A$}
 	\RightLabel{EquivConst}
 	\UnaryInfC{$A\leftrightarrow B$}
 	\DisplayProof
 	\hskip 1.5em
 	\AxiomC{$A \leftrightarrow B$}
 	\RightLabel{EquivExp}
 	\UnaryInfC{$A\rightarrow B \enspace\&\enspace B \rightarrow A$}
 	\DisplayProof
 	
 \end{center}
 
 \begin{remark}
 	Certain combinations of rules occur frequently and it is convenient to have dervived rules(or shortcuts) such as
 	
 	\vspace{5mm}
 	
 	\AxiomC{$A \rightarrow B$}
 	\AxiomC{$B\rightarrow A$}
 	\RightLabel{EquivJoin}
 	\BinaryInfC{$A\leftrightarrow B$}
 	\DisplayProof
 	
 	\vspace{5mm}

 	for
 	
 	\vspace{5mm}
 	
 	\AxiomC{$A \rightarrow B$}
 	\AxiomC{$B\rightarrow A$}
 	\RightLabel{AndInt}
 	\BinaryInfC{$A\rightarrow B \enspace\&\enspace B \rightarrow A$}
 	\RightLabel{EquivCont}
 	\UnaryInfC{$A\leftrightarrow B$}
 	\DisplayProof
 	\vspace{5mm}
 	
 	and two other obvious shortcuts EquivRight and EquivLeft.
 Also the important derived rule

 	\vspace{5mm}
 	
 	\AxiomC{$A$}
 	\RightLabel{FreeSub$_{y,t}$}
 	\UnaryInfC{$A[t/y]$}
 	\DisplayProof
 	
 	\vspace{5mm}
 	
 	for
 	\vspace{5mm}
 	
 	\AxiomC{$A$}
 	\RightLabel{ForallInt$_y$}
 	\UnaryInfC{$\forall x. A[x/y]$}
 	\RightLabel{ForallElim}
 	\UnaryInfC{$A[t/x]$}
 	\DisplayProof

 \end{remark}

\section*{Linearised Natural Deduction}
A linear proof in PyLog (for a given proof environment) is a list of \emph{proof elements}. Each proof element $p$ is a triple $(A,par,dis)$ where $A$ is a formula and $par$ and $dis$ are lists of integers. If $p$ occurs in position $m$  (we say that $m$ is $p$'s number) and the list has length $m$ then the elements of $par$ (parents) must be strictly less than $n$ and those $dis$(discharges) strictly larger than $n$. Rules are applied by adding a new proof element to the end of the list and possibly updating previous proof entries.
Given a linear proof and an element $p$ we can recursively backtrack the parents to obtain a \emph{dependency tree} of proof-elements. The dependencies of $p$ are obtained by taking the set of leaves of this tree are removing those proof-elements having $dis$ containing a number which occurs in the dependency tree. It is also easy to see how given a proof tree we can obtain a linear proof.
In PyLog rules have the general format

\[ (Name, Parents, Parameters) \]

where Name is the name of the rule, Parents is the list of the number previous elements which the rule is applied to and Parameters can constain formulas, terms, variables, position lists, etc.
In PyLog rules are entered as  Python function. The arguments will specify all the required information in Parents and Parameters.

The PyLog command Qed(ForNum) checks if the formula has discharged all its assumptions.

It is helpful to look at a snippet in the definition of some classes:

\begin{verbatim}
class ProofElement:
 def __init__(self,name,dependencies,parameters,discharging,formula):
  self.name = name
  self.dependencies = dependencies
  self.parameters = parameters
  self.discharging = discharging
  self.formula = formula
  self.dischargedby = []
  self.pos = 0
  self.qed=False 
  self.comment=""

class ProofEnvironment:
 def __init__(self,proof,name):
  self.proof = proof
  self.name = name
  self.definitions = {}
  self.definitionequations = []
  self.axioms = []
  self.theorems = []
  self.log = []

 def CheckRange(self, dependencies):
  for dep in dependencies:
   if dep > len(self.proof):
    return False
   return True

 def GetTree(self,proofelement):
  out = [proofelement.pos-1]
  for dep in proofelement.dependencies:
    out = out + self.GetTree(self.proof[dep])
    return out

 def GetHyp(self,proofelement):
   if proofelement.name =="Hyp":
    return [proofelement.pos-1]
   out = []
   for dep in proofelement.dependencies:
     out = out + self.GetHyp(self.proof[dep])
   return out

 def CheckDischargedBy(self, hyp, proofelem):
   if len(self.proof[hyp].dischargedby) ==0:
     return False
   for h-1 in self.proof[hyp].dischargedby:
     if h in self.GetTree(self.proof[proofelem]):
       return True
     return False

 def GetHypDep(self,proofelement): 
   aux = []
   for  h in self.GetHyp(proofelement):
     if len(Intersect([x-1 for x in self.proof[h].dischargedby],
      self.GetTree(proofelement))) == 0:
       aux.append(h)
   return aux 

   (...)
\end{verbatim}	

\section*{List of Core PyLog Rules}

\begin{verbatim}
Logical Rules
=============

AndInt(ForNum, ForNum)

AndElimL(ForNum)

AndElimR(ForNum)

ImpInt(ForNum, DisNum)

ImpElim(ForNum,ImpNum)

OrIntL(ForNum,formula)

OrIntR(ForNum,formula)

OrElim(OrForNum, LeftHypNum, LeftConNum, RightHypNum, RightConNum)

ForallInt(ForNum, VarName, newVarName)

ForallElim(ForNum, term)

ExistsInt(ForNum, term, newVarName, PositionList)

ExistsElim(ExistsForNum, InstForNum,ConForNum, instVariable)

AbsI(BotForNum)

AbsC(NegForNum, BotConNum)


Class Rules
=========

ClassElim(MemForNum)

ClassInt(ForNum, newVarName)



Equality Rules
===============

Identity(term)

Symmetry(EqForNum)

EqualitySub(ForNum, EqForNum, PositionList)





Second-Order Rule
=================

PolySub(ForNum, SecondOrderVarName, formula)


Proof Environment Rules
=======================

AxInt(number)

TheoremInt(number)

DefEqInt(number)

DefExp(ForNum, predicateName, PositionList)

DefSub(ForNum, predicateName, ArgList, PositionList)


Other Rules
==========

Qed(ForNum)

EquivConst(ForNum)

EquivExp(ForNum)

EquivLeft(ForNum)

EquivRight(ForNum)

FreeSub(ForNum, VarName, Term)

\end{verbatim}

\section*{Pylog Commands}

We have a list of commands for setting up the proof environment, that is, for introducing the axioms, assumed theorems, defined constants and symbols and defined predicates that will be used in the proof.

\begin{verbatim}
	
Hyp(Formula)	
	
NewAx(Formula)

AddPredicate(PredicateName, Arity, PrefixBook)

NewDef(PredicateName, ArgList, Formula)

AddConstants(NameList)

AddFunction(FunctionName, Arity, PrefixBool)

NewDefEq(EquationFormula)

AddTheorem(Formula)

\end{verbatim}

Hyp introduces a formula as a hypothesis.

When defining functions with NewDefEq() we must first use AddFunction() specifying the name, arity and whether the function is to be displayed with prefix or infix notation (for binary functions). The argument for NewDefEq must not have the exterior parenthesis. For instance

\begin{verbatim}

NewDefEq("rus = extension z. neg Elem(z,z)")

\end{verbatim}

 For constants we use AddConstant().
We also must take care that we have enough variables via the AddVariables(VarList) function.

Then we have a list of commands which displays information about the current proof and proof environment. ShowDefinitions()  displays the defined predicates. ShowDefEquations() displays the defined constants and functions. ShowAxioms() displays the axioms.  ShowTheorems() displays the assumed theorems which may be used in the proof (it is not advisable to alter this list during the proof). ShowProof() displays the current state of the proof and ShowLog() shows the list of previous succesful rule commands which constitute the proof. We also have a command Undo() which deletes the last element of the proof. Hypotheses(n) shows the hypotheses which formula $n$ depends on.

If a theorem has already been saved you can view the conclusion with ViewTheorem(Name) or add it directly to the proof environment with the LoadTheorem(Name) command - \emph{provided that the required environment has been previously loaded}.

\section*{Using Pylog}

In PyLog a \emph{theory} is a directory whose files are \emph{theorems}. A theorem consists of both a proof environment and a proof in this environment (either complete or incomplete).  All theorems in a theory should ideally have the same proof environment. The theorem to be proved ideally should occur at the end of the proof and have been tested with the command Qed(Number).
There should also be an "empty" theorem which is to be seen as the proof environment that must be loaded in order to start writing a new theorem. The command Load(Name) loads a proof environment or theorem and the command Save(Name) will save the current proof environment or theorem. The command ViewTheorem(Name) will not load anything but only display the last line of the proof. The command ViewTheory(DirName) will likewise display all the theorems in the directory.

To use Pylog Python 3.$*$ is required. PyLog runs from a terminal through the Python CLI.
  Clone the repository\footnote{\url{https://github.com/owl77/PyLog}} on GitHub, enter the folder, and enter
    \begin{verbatim}
$ python -i proofenvironment.py

Welcome to PyLog 1.0

Natural Deduction Proof Assistant and Proof Checker

(c) 2020  C. Lewis Protin

>>> 
 \end{verbatim}

Our project is to have a complete verified formalisation of all the theorems of Set Theory in the Appendix of \cite{kel}.

In the PyLog folder we have the saved Kelley-Morse environment. We load this by Load("Kelley-Morse"). When a command is succesful PyLog will return True. We can now examine the axioms and definitions:

\begin{verbatim}
>>> ShowAxioms()
0. ∀x.∀y.((x = y) <-> ∀z.((z ε x) <-> (z ε y))) 
1. Set(x) -> ∃y.(Set(y) & ∀z.((z ⊂ x) -> (z ε y))) 
2. (Set(x) & Set(y)) -> Set((x ∪ y)) 
3. (Function(f) & Set(domain(f))) -> Set(range(f)) 
4. Set(x) -> Set(∪x) 
5. ¬(x = 0) -> ∃y.((y ε x) & ((y ∩ x) = 0)) 
6. ∃y.((Set(y) & (0 ε y)) & ∀x.((x ε y) -> (suc x ε y))) 
7. ∃f.(Choice(f) & (domain(f) = (U ~ {0}))) 
>>> ShowDefEquations()
0. (x ∪ y) = {z: ((z ε x) v (z ε y))} 
1. (x ∩ y) = {z: ((z ε x) & (z ε y))} 
2. ~x = {y: ¬(y ε x)} 
3. (x ~ y) = (x ∩ ~y) 
4. 0 = {x: ¬(x = x)} 
5. U = {x: (x = x)} 
6. ∪x = {z: ∃y.((y ε x) & (z ε y))} 
7. ∩x = {z: ∀y.((y ε x) -> (z ε y))} 
8. Px = {y: (y ⊂ x)} 
9. {x} = {z: ((z ε U) -> (z = x))} 
10. {x,y} = ({x} ∪ {y}) 
11. (x,y) = {x,{x,y}} 
12. proj1(x) = ∩∩x 
13. proj2(x) = (∩∪x ∪ (∪∪x ~ ∪∩x)) 
14. (a∘b) = {w: ∃x.∃y.∃z.((((x,y) ε a) & ((y,z) ε b)) & (w = (x,z)))} 
15. (r)⁻¹ = {z: ∃x.∃y.(((x,y) ε r) & (z = (y,x)))} 
16. domain(f) = {x: ∃y.((x,y) ε f)} 
17. range(f) = {y: ∃x.((x,y) ε f)} 
18. (f'x) = ∩{y: ((x,y) ε f)} 
19. (x X y) = {z: ∃a.∃b.((z = (a,b)) & ((a ε x) & (b ε y)))} 
20. func(x,y) = {f: (Function(f) & ((domain(f) = x) & (range(f) = y)))} 
21. E = {z: ∃x.∃y.((z = (x,y)) & (x ε y))} 
22. ord = {x: Ordinal(x)} 
23. suc x = (x ∪ {x}) 
24. (f|x) = (f ∩ (x X U)) 
25. ω = {x: Integer(x)} 
>>> ShowDefinitions()
Set(x) <-> ∃y.(x ε y)
(x ⊂ y) <-> ∀z.((z ε x) -> (z ε y))
Relation(r) <-> ∀z.((z ε r) -> ∃x.∃y.(z = (x,y)))
Function(f) <-> (Relation(f) & ∀x.∀y.∀z.((((x,y) ε f) & ((x,z) ε f)) -> (y = z)))
Trans(r) <-> ∀x.∀y.∀z.((((x,y) ε r) & ((y,z) ε r)) -> ((x,z) ε r))
Connects(r,x) <-> ∀y.∀z.(((y ε x) & (z ε x)) -> ((y = z) v (((y,z) ε r) v ((z,y) ε r))))
Asymmetric(r,x) <-> ∀y.∀z.(((y ε x) & (z ε x)) -> (((y,z) ε r) -> ¬((z,y) ε r)))
First(r,x,z) <-> ((z ε x) & ∀y.((y ε x) -> ¬((y,z) ε r)))
WellOrders(r,x) <-> (Connects(r,x) & ∀y.(((y ⊂ x) & ¬(y = 0)) -> ∃z.First(r,y,z)))
Section(r,x,y) <-> (((y ⊂ x) & WellOrders(r,x)) & ∀u.∀v.((((u ε x)
 & (v ε y)) & ((u,v) ε r)) -> (u ε y)))
OrderPreserving(f,r,s) <-> ((Function(f) & (WellOrders(r,domain(f))
 & WellOrders(r,range(f)))) & ∀u.∀v.((((u ε domain(f)) 
 & (v ε domain(f))) & ((u,v) ε r)) -> (((f'u),(f'v)) ε r)))
1-to-1(f) <-> (Function(f) & Function((f)⁻¹))
Full(x) <-> ∀y.((y ε x) -> (y ⊂ x))
Ordinal(x) <-> (Full(x) & Connects(E,x))
Integer(x) <-> (Ordinal(x) & WellOrders((E)⁻¹,x))
Choice(f) <-> (Function(f) & ∀y.((y ε domain(f)) -> ((f'y) ε y)))
Equi(x,y) <-> ∃f.(1-to-1(f) & ((domain(f) = x) & (range(f) = y)))
Card(x) <-> (Ordinal(x) & ∀y.(((y ε x) & (y ε ord)) -> ¬Equi(y,x)))
TransIn(r,x) <-> ∀u.∀v.∀w.(((u ε x) & ((v ε x) & (w ε x))) ->
 ((((u,v) ε r) & ((v,w) ε r)) -> ((u,w) ε r)))
\end{verbatim}

We can also check by ShowProof() that the proof is empty. By default expressions are displayed using pretty printing (Unicode character) which can use infix notation.  The pretty printing can be changed via parser.prettyprint[FunctionNameString] = PrettyString. Expressions are entered in a strictly functional way (with the exception of logical connectives, extensions and quantifiers). 

\begin{verbatim}
Input and default "pretty" display
===================================
neg A                     ¬A  
bigunion(x)               ∪x
bigintersection           ∩x
union(x,y)               (x ∪ y)
intersection(x,y)        (x ∩ y)
extension x. A           {x: A}
forall x. A              ∀x. A
exists x. A              ∃x. A
Elem(x,y)               (x ε y)
app(f,x)                 (f'x)
pair(x,y)                {x,y}
singleton(x)             {x}
orderedpair(x,y)         (x,y)
prod(x,y)                (x X y)
complement1(x)              ~x
complement2(x)            (x ~ y)
parts(x)                  Px
comp(a,b)                 (a∘b) 
inv(r)                     (r)⁻¹ 
restrict(f,x)              (f|x)
int                         ω

\end{verbatim}

Note that $\neg A$ is the pretty print for $A\rightarrow \bot$. When using the rules of Pylog
we must think of $\neg A$ this way.
Conjunction $\&$ is usually entered in infix style $(A\enspace \& \enspace B)$ but PyLog will create and group parenthesis to the right: thus $(A \enspace \& \enspace B\enspace  \& \enspace C)$ is interpreted as $(A\enspace \& \enspace (B \enspace \& \enspace C))$.

\section*{First Proof in Pylog}

In this section we prove theorem 4 of \cite{kel} in the Kelley-Morse proof environment:
\begin{verbatim}
((z ε (x ∪ y)) <-> ((z ε x) v (z ε y))) & ((z ε (x ∩ y)) <-> ((z ε x) & (z ε y)))
\end{verbatim}

We give here full  details of a session in which we prove the first half of the theorem.

\begin{verbatim}
Welcome to PyLog 1.0

Natural Deduction Proof Assistant and Proof Checker

(c) 2020  C. Lewis Protin

>>> Load("Kelley-Morse")
True
>>> ShowProof()
>>> Hyp("Elem(z,union(x,y))")
0. z ε (x ∪ y) Hyp 
True
>>> ShowDefEquations()
0. (x ∪ y) = {z: ((z ε x) v (z ε y))} 
1. (x ∩ y) = {z: ((z ε x) & (z ε y))} 
2. ~x = {y: ¬(y ε x)} 
3. (x ~ y) = (x ∩ ~y) 
4. 0 = {x: ¬(x = x)} 
5. U = {x: (x = x)} 
6. ∪x = {z: ∃y.((y ε x) & (z ε y))} 
7. ∩x = {z: ∀y.((y ε x) -> (z ε y))} 
8. Px = {y: (y ⊂ x)} 
9. {x} = {z: ((z ε U) -> (z = x))} 
10. {x,y} = ({x} ∪ {y}) 
11. (x,y) = {x,{x,y}} 
12. proj1(x) = ∩∩x 
13. proj2(x) = (∩∪x ∪ (∪∪x ~ ∪∩x)) 
14. (a∘b) = {w: ∃x.∃y.∃z.((((x,y) ε a) & ((y,z) ε b)) & (w = (x,z)))} 
15. (r)⁻¹ = {z: ∃x.∃y.(((x,y) ε r) & (z = (y,x)))} 
16. domain(f) = {x: ∃y.((x,y) ε f)} 
17. range(f) = {y: ∃x.((x,y) ε f)} 
18. (f'x) = ∩{y: ((x,y) ε f)} 
19. (x X y) = {z: ∃a.∃b.((z = (a,b)) & ((a ε x) & (b ε y)))} 
20. func(x,y) = {f: (Function(f) & ((domain(f) = x) & (range(f) = y)))} 
21. E = {z: ∃x.∃y.((z = (x,y)) & (x ε y))} 
22. ord = {x: Ordinal(x)} 
23. suc x = (x ∪ {x}) 
24. (f|x) = (f ∩ (x X U)) 
25. ω = {x: Integer(x)} 
>>> DefEqInt(0)
0. z ε (x ∪ y) Hyp 
1. (x ∪ y) = {z: ((z ε x) v (z ε y))} DefEqInt 
True
>>> EqualitySub(0,1,[0])
0. z ε (x ∪ y) Hyp 
1. (x ∪ y) = {z: ((z ε x) v (z ε y))} DefEqInt 
2. z ε {z: ((z ε x) v (z ε y))} EqualitySub 0 1
True
>>> ClassElim(2)
0. z ε (x ∪ y) Hyp 
1. (x ∪ y) = {z: ((z ε x) v (z ε y))} DefEqInt 
2. z ε {z: ((z ε x) v (z ε y))} EqualitySub 0 1
3. Set(z) & ((z ε x) v (z ε y)) ClassElim 2
True
>>> AndElimR(3)
0. z ε (x ∪ y) Hyp 
1. (x ∪ y) = {z: ((z ε x) v (z ε y))} DefEqInt 
2. z ε {z: ((z ε x) v (z ε y))} EqualitySub 0 1
3. Set(z) & ((z ε x) v (z ε y)) ClassElim 2
4. (z ε x) v (z ε y) AndElimR 3
True
>>> ImpInt(4,0)
0. z ε (x ∪ y) Hyp 
1. (x ∪ y) = {z: ((z ε x) v (z ε y))} DefEqInt 
2. z ε {z: ((z ε x) v (z ε y))} EqualitySub 0 1
3. Set(z) & ((z ε x) v (z ε y)) ClassElim 2
4. (z ε x) v (z ε y) AndElimR 3
5. (z ε (x ∪ y)) -> ((z ε x) v (z ε y)) ImpInt 4
True
>>> Qed(5)
(...)
5. (z ε (x ∪ y)) -> ((z ε x) v (z ε y)) ImpInt 4 Qed
>>> Hyp("(Elem(z,x) v Elem(z,y))")
0. z ε (x ∪ y) Hyp 
1. (x ∪ y) = {z: ((z ε x) v (z ε y))} DefEqInt 
2. z ε {z: ((z ε x) v (z ε y))} EqualitySub 0 1
3. Set(z) & ((z ε x) v (z ε y)) ClassElim 2
4. (z ε x) v (z ε y) AndElimR 3
5. (z ε (x ∪ y)) -> ((z ε x) v (z ε y)) ImpInt 4 Qed
6. (z ε x) v (z ε y) Hyp 
True
>>> Hyp("Elem(z,x)")
0. z ε (x ∪ y) Hyp 
1. (x ∪ y) = {z: ((z ε x) v (z ε y))} DefEqInt 
2. z ε {z: ((z ε x) v (z ε y))} EqualitySub 0 1
3. Set(z) & ((z ε x) v (z ε y)) ClassElim 2
4. (z ε x) v (z ε y) AndElimR 3
5. (z ε (x ∪ y)) -> ((z ε x) v (z ε y)) ImpInt 4 Qed
6. (z ε x) v (z ε y) Hyp 
7. z ε x Hyp 
True
>>> ExistsInt(7,"x","x",[0])
0. z ε (x ∪ y) Hyp 
1. (x ∪ y) = {z: ((z ε x) v (z ε y))} DefEqInt 
2. z ε {z: ((z ε x) v (z ε y))} EqualitySub 0 1
3. Set(z) & ((z ε x) v (z ε y)) ClassElim 2
4. (z ε x) v (z ε y) AndElimR 3
5. (z ε (x ∪ y)) -> ((z ε x) v (z ε y)) ImpInt 4 Qed
6. (z ε x) v (z ε y) Hyp 
7. z ε x Hyp 
8. ∃x.(z ε x) ExistsInt 7
True
>>> ShowDefinitions()
Set(x) <-> ∃y.(x ε y)
(x ⊂ y) <-> ∀z.((z ε x) -> (z ε y))
Relation(r) <-> ∀z.((z ε r) -> ∃x.∃y.(z = (x,y)))
Function(f) <-> (Relation(f) & ∀x.∀y.∀z.((((x,y) ε f) & ((x,z) ε f)) -> (y = z)))
Trans(r) <-> ∀x.∀y.∀z.((((x,y) ε r) & ((y,z) ε r)) -> ((x,z) ε r))
Connects(r,x) <-> ∀y.∀z.(((y ε x) & (z ε x)) -> ((y = z) v (((y,z) ε r) v ((z,y) ε r))))
Asymmetric(r,x) <-> ∀y.∀z.(((y ε x) & (z ε x)) -> (((y,z) ε r) -> ¬((z,y) ε r)))
First(r,x,z) <-> ((z ε x) & ∀y.((y ε x) -> ¬((y,z) ε r)))
WellOrders(r,x) <-> (Connects(r,x) & ∀y.(((y ⊂ x) & ¬(y = 0)) -> ∃z.First(r,y,z)))
Section(r,x,y) <-> (((y ⊂ x) & WellOrders(r,x)) & ∀u.∀v.((((u ε x) 
& (v ε y)) & ((u,v) ε r)) -> (u ε y)))
OrderPreserving(f,r,s) <-> ((Function(f) & (WellOrders(r,domain(f))
 & WellOrders(r,range(f)))) & ∀u.∀v.((((u ε domain(f)) & (v ε domain(f))) & ((u,v) ε r)) -> (((f'u),(f'v)) ε r)))
1-to-1(f) <-> (Function(f) & Function((f)⁻¹))
Full(x) <-> ∀y.((y ε x) -> (y ⊂ x))
Ordinal(x) <-> (Full(x) & Connects(E,x))
Integer(x) <-> (Ordinal(x) & WellOrders((E)⁻¹,x))
Choice(f) <-> (Function(f) & ∀y.((y ε domain(f)) -> ((f'y) ε y)))
Equi(x,y) <-> ∃f.(1-to-1(f) & ((domain(f) = x) & (range(f) = y)))
Card(x) <-> (Ordinal(x) & ∀y.(((y ε x) & (y ε ord)) -> ¬Equi(y,x)))
TransIn(r,x) <-> ∀u.∀v.∀w.(((u ε x) & ((v ε x) & (w ε x))) -> ((((u,v) ε r) 
& ((v,w) ε r)) -> ((u,w) ε r)))
>>> DefSub(8,"Set","z",[0])
0. z ε (x ∪ y) Hyp 
1. (x ∪ y) = {z: ((z ε x) v (z ε y))} DefEqInt 
2. z ε {z: ((z ε x) v (z ε y))} EqualitySub 0 1
3. Set(z) & ((z ε x) v (z ε y)) ClassElim 2
4. (z ε x) v (z ε y) AndElimR 3
5. (z ε (x ∪ y)) -> ((z ε x) v (z ε y)) ImpInt 4 Qed
6. (z ε x) v (z ε y) Hyp 
7. z ε x Hyp 
8. ∃x.(z ε x) ExistsInt 7
9. Set(z) DefSub 8
True
>>> Hyp("Elem(z,y)")
0. z ε (x ∪ y) Hyp 
1. (x ∪ y) = {z: ((z ε x) v (z ε y))} DefEqInt 
2. z ε {z: ((z ε x) v (z ε y))} EqualitySub 0 1
3. Set(z) & ((z ε x) v (z ε y)) ClassElim 2
4. (z ε x) v (z ε y) AndElimR 3
5. (z ε (x ∪ y)) -> ((z ε x) v (z ε y)) ImpInt 4 Qed
6. (z ε x) v (z ε y) Hyp 
7. z ε x Hyp 
8. ∃x.(z ε x) ExistsInt 7
9. Set(z) DefSub 8
10. z ε y Hyp 
True
>>> ExistsInt(10, "y","x",[0])
0. z ε (x ∪ y) Hyp 
1. (x ∪ y) = {z: ((z ε x) v (z ε y))} DefEqInt 
2. z ε {z: ((z ε x) v (z ε y))} EqualitySub 0 1
3. Set(z) & ((z ε x) v (z ε y)) ClassElim 2
4. (z ε x) v (z ε y) AndElimR 3
5. (z ε (x ∪ y)) -> ((z ε x) v (z ε y)) ImpInt 4 Qed
6. (z ε x) v (z ε y) Hyp 
7. z ε x Hyp 
8. ∃x.(z ε x) ExistsInt 7
9. Set(z) DefSub 8
10. z ε y Hyp 
11. ∃x.(z ε x) ExistsInt 10
True
>>> DefSub(11,"Set","z",[0])
0. z ε (x ∪ y) Hyp 
1. (x ∪ y) = {z: ((z ε x) v (z ε y))} DefEqInt 
2. z ε {z: ((z ε x) v (z ε y))} EqualitySub 0 1
3. Set(z) & ((z ε x) v (z ε y)) ClassElim 2
4. (z ε x) v (z ε y) AndElimR 3
5. (z ε (x ∪ y)) -> ((z ε x) v (z ε y)) ImpInt 4 Qed
6. (z ε x) v (z ε y) Hyp 
7. z ε x Hyp 
8. ∃x.(z ε x) ExistsInt 7
9. Set(z) DefSub 8
10. z ε y Hyp 
11. ∃x.(z ε x) ExistsInt 10
12. Set(z) DefSub 11
True
>>> OrElim(6,7,9,10,12)
0. z ε (x ∪ y) Hyp 
1. (x ∪ y) = {z: ((z ε x) v (z ε y))} DefEqInt 
2. z ε {z: ((z ε x) v (z ε y))} EqualitySub 0 1
3. Set(z) & ((z ε x) v (z ε y)) ClassElim 2
4. (z ε x) v (z ε y) AndElimR 3
5. (z ε (x ∪ y)) -> ((z ε x) v (z ε y)) ImpInt 4 Qed
6. (z ε x) v (z ε y) Hyp 
7. z ε x Hyp 
8. ∃x.(z ε x) ExistsInt 7
9. Set(z) DefSub 8
10. z ε y Hyp 
11. ∃x.(z ε x) ExistsInt 10
12. Set(z) DefSub 11
13. Set(z) OrElim 6 7 9 10 12
True
>>> AndInt(13,6)
0. z ε (x ∪ y) Hyp 
1. (x ∪ y) = {z: ((z ε x) v (z ε y))} DefEqInt 
2. z ε {z: ((z ε x) v (z ε y))} EqualitySub 0 1
3. Set(z) & ((z ε x) v (z ε y)) ClassElim 2
4. (z ε x) v (z ε y) AndElimR 3
5. (z ε (x ∪ y)) -> ((z ε x) v (z ε y)) ImpInt 4 Qed
6. (z ε x) v (z ε y) Hyp 
7. z ε x Hyp 
8. ∃x.(z ε x) ExistsInt 7
9. Set(z) DefSub 8
10. z ε y Hyp 
11. ∃x.(z ε x) ExistsInt 10
12. Set(z) DefSub 11
13. Set(z) OrElim 6 7 9 10 12
14. Set(z) & ((z ε x) v (z ε y)) AndInt 13 6
True
>>> ClassInt(14,"z")
0. z ε (x ∪ y) Hyp 
1. (x ∪ y) = {z: ((z ε x) v (z ε y))} DefEqInt 
2. z ε {z: ((z ε x) v (z ε y))} EqualitySub 0 1
3. Set(z) & ((z ε x) v (z ε y)) ClassElim 2
4. (z ε x) v (z ε y) AndElimR 3
5. (z ε (x ∪ y)) -> ((z ε x) v (z ε y)) ImpInt 4 Qed
6. (z ε x) v (z ε y) Hyp 
7. z ε x Hyp 
8. ∃x.(z ε x) ExistsInt 7
9. Set(z) DefSub 8
10. z ε y Hyp 
11. ∃x.(z ε x) ExistsInt 10
12. Set(z) DefSub 11
13. Set(z) OrElim 6 7 9 10 12
14. Set(z) & ((z ε x) v (z ε y)) AndInt 13 6
15. z ε {z: ((z ε x) v (z ε y))} ClassInt 14
True
>>> Symmetry(1)
0. z ε (x ∪ y) Hyp 
1. (x ∪ y) = {z: ((z ε x) v (z ε y))} DefEqInt 
2. z ε {z: ((z ε x) v (z ε y))} EqualitySub 0 1
3. Set(z) & ((z ε x) v (z ε y)) ClassElim 2
4. (z ε x) v (z ε y) AndElimR 3
5. (z ε (x ∪ y)) -> ((z ε x) v (z ε y)) ImpInt 4 Qed
6. (z ε x) v (z ε y) Hyp 
7. z ε x Hyp 
8. ∃x.(z ε x) ExistsInt 7
9. Set(z) DefSub 8
10. z ε y Hyp 
11. ∃x.(z ε x) ExistsInt 10
12. Set(z) DefSub 11
13. Set(z) OrElim 6 7 9 10 12
14. Set(z) & ((z ε x) v (z ε y)) AndInt 13 6
15. z ε {z: ((z ε x) v (z ε y))} ClassInt 14
16. {z: ((z ε x) v (z ε y))} = (x ∪ y) Symmetry 1
True
>>> EqualitySub(15,16,[0])
0. z ε (x ∪ y) Hyp 
1. (x ∪ y) = {z: ((z ε x) v (z ε y))} DefEqInt 
2. z ε {z: ((z ε x) v (z ε y))} EqualitySub 0 1
3. Set(z) & ((z ε x) v (z ε y)) ClassElim 2
4. (z ε x) v (z ε y) AndElimR 3
5. (z ε (x ∪ y)) -> ((z ε x) v (z ε y)) ImpInt 4 Qed
6. (z ε x) v (z ε y) Hyp 
7. z ε x Hyp 
8. ∃x.(z ε x) ExistsInt 7
9. Set(z) DefSub 8
10. z ε y Hyp 
11. ∃x.(z ε x) ExistsInt 10
12. Set(z) DefSub 11
13. Set(z) OrElim 6 7 9 10 12
14. Set(z) & ((z ε x) v (z ε y)) AndInt 13 6
15. z ε {z: ((z ε x) v (z ε y))} ClassInt 14
16. {z: ((z ε x) v (z ε y))} = (x ∪ y) Symmetry 1
17. z ε (x ∪ y) EqualitySub 15 16
True
>>> ImpInt(17,6)
0. z ε (x ∪ y) Hyp 
1. (x ∪ y) = {z: ((z ε x) v (z ε y))} DefEqInt 
2. z ε {z: ((z ε x) v (z ε y))} EqualitySub 0 1
3. Set(z) & ((z ε x) v (z ε y)) ClassElim 2
4. (z ε x) v (z ε y) AndElimR 3
5. (z ε (x ∪ y)) -> ((z ε x) v (z ε y)) ImpInt 4 Qed
6. (z ε x) v (z ε y) Hyp 
7. z ε x Hyp 
8. ∃x.(z ε x) ExistsInt 7
9. Set(z) DefSub 8
10. z ε y Hyp 
11. ∃x.(z ε x) ExistsInt 10
12. Set(z) DefSub 11
13. Set(z) OrElim 6 7 9 10 12
14. Set(z) & ((z ε x) v (z ε y)) AndInt 13 6
15. z ε {z: ((z ε x) v (z ε y))} ClassInt 14
16. {z: ((z ε x) v (z ε y))} = (x ∪ y) Symmetry 1
17. z ε (x ∪ y) EqualitySub 15 16
18. ((z ε x) v (z ε y)) -> (z ε (x ∪ y)) ImpInt 17
True
>>> Qed(18)
True
>>> ShowProof()
(...)
18. ((z ε x) v (z ε y)) -> (z ε (x ∪ y)) ImpInt 17 Qed
>>> AndInt(5,18)
0. z ε (x ∪ y) Hyp 
1. (x ∪ y) = {z: ((z ε x) v (z ε y))} DefEqInt 
2. z ε {z: ((z ε x) v (z ε y))} EqualitySub 0 1
3. Set(z) & ((z ε x) v (z ε y)) ClassElim 2
4. (z ε x) v (z ε y) AndElimR 3
5. (z ε (x ∪ y)) -> ((z ε x) v (z ε y)) ImpInt 4 Qed
6. (z ε x) v (z ε y) Hyp 
7. z ε x Hyp 
8. ∃x.(z ε x) ExistsInt 7
9. Set(z) DefSub 8
10. z ε y Hyp 
11. ∃x.(z ε x) ExistsInt 10
12. Set(z) DefSub 11
13. Set(z) OrElim 6 7 9 10 12
14. Set(z) & ((z ε x) v (z ε y)) AndInt 13 6
15. z ε {z: ((z ε x) v (z ε y))} ClassInt 14
16. {z: ((z ε x) v (z ε y))} = (x ∪ y) Symmetry 1
17. z ε (x ∪ y) EqualitySub 15 16
18. ((z ε x) v (z ε y)) -> (z ε (x ∪ y)) ImpInt 17 Qed
19. ((z ε (x ∪ y)) -> ((z ε x) v (z ε y))) & 
(((z ε x) v (z ε y)) -> (z ε (x ∪ y))) AndInt 5 18
True
>>> EquivConst(19)
0. z ε (x ∪ y) Hyp 
1. (x ∪ y) = {z: ((z ε x) v (z ε y))} DefEqInt 
2. z ε {z: ((z ε x) v (z ε y))} EqualitySub 0 1
3. Set(z) & ((z ε x) v (z ε y)) ClassElim 2
4. (z ε x) v (z ε y) AndElimR 3
5. (z ε (x ∪ y)) -> ((z ε x) v (z ε y)) ImpInt 4 Qed
6. (z ε x) v (z ε y) Hyp 
7. z ε x Hyp 
8. ∃x.(z ε x) ExistsInt 7
9. Set(z) DefSub 8
10. z ε y Hyp 
11. ∃x.(z ε x) ExistsInt 10
12. Set(z) DefSub 11
13. Set(z) OrElim 6 7 9 10 12
14. Set(z) & ((z ε x) v (z ε y)) AndInt 13 6
15. z ε {z: ((z ε x) v (z ε y))} ClassInt 14
16. {z: ((z ε x) v (z ε y))} = (x ∪ y) Symmetry 1
17. z ε (x ∪ y) EqualitySub 15 16
18. ((z ε x) v (z ε y)) -> (z ε (x ∪ y)) ImpInt 17 Qed
19. ((z ε (x ∪ y)) -> ((z ε x) v (z ε y))) &
 (((z ε x) v (z ε y)) -> (z ε (x ∪ y))) AndInt 5 18
20. (z ε (x ∪ y)) <-> ((z ε x) v (z ε y)) EquivConst 
True
>>> Qed(20)
(...)
20. (z ε (x ∪ y)) <-> ((z ε x) v (z ε y)) EquivConst  Qed
>>> Save("Th4")
True
\end{verbatim}

The proof is fully codified by the log (and the proof environment):

\begin{verbatim}
>>> ShowLog()
0. Hyp("Elem(z, union(x,y))")
1. DefEqInt(0)
2. EqualitySub(0,1,[0])
3. ClassElim(2)
4. AndElimR(3)
5. ImpInt(4,0)
6. Hyp("(Elem(z,x) v Elem(z,y))")
7. Hyp("Elem(z,x)")
8. ExistsInt(7,"x","x",[0])
9. DefSub(8,"Set",["z"],[0])
10. Hyp("Elem(z,y)")
11. ExistsInt(10,"y","y",[0])
12. DefSub(11,"Set",["z"],[0])
13. OrElim(6,7,9,10,12)
14. AndInt(13,6)
15. ClassInt(14,"z")
16. Symmetry(1)
17. EqualitySub(15,16,[0])
18. ImpInt(17,6)
19. AndInt(5,18)
20. EquivConst(19)
\end{verbatim}	

The command GenerateProof() deletes the proof and then generates it again from the Proof.log  field and runs Qed() on the last line. This allows
any proof log to be formally checked for a given environement. The log  contains all the information needed to generate the proof.
The command UsedTheorems() gives the list of other theorems used in the proof.

Similarly the second half of the conjunction is proven:

\begin{verbatim}
21. z ε (x ∩ y) Hyp 
22. (x ∩ y) = {z: ((z ε x) & (z ε y))} DefEqInt 
23. z ε {z: ((z ε x) & (z ε y))} EqualitySub 21 22
24. Set(z) & ((z ε x) & (z ε y)) ClassElim 23
25. (z ε x) & (z ε y) AndElimR 24
26. (z ε (x ∩ y)) -> ((z ε x) & (z ε y)) ImpInt 25 Qed
27. (z ε x) & (z ε y) Hyp 
28. z ε x AndElimL 27
29. ∃x.(z ε x) ExistsInt 28
30. Set(z) DefSub 29
31. Set(z) & ((z ε x) & (z ε y)) AndInt 30 27
32. z ε {z: ((z ε x) & (z ε y))} ClassInt 31
33. {z: ((z ε x) & (z ε y))} = (x ∩ y) Symmetry 22
34. z ε (x ∩ y) EqualitySub 32 33
35. ((z ε x) & (z ε y)) -> (z ε (x ∩ y)) ImpInt 34 Qed
36. ((z ε (x ∩ y)) -> ((z ε x) & (z ε y))) & (((z ε x) & (z ε y)) -> (z ε (x ∩ y)))
 AndInt 26 35
37. (z ε (x ∩ y)) <-> ((z ε x) & (z ε y)) EquivConst  Qed
38. ((z ε (x ∪ y)) <-> ((z ε x) v (z ε y))) & ((z ε (x ∩ y)) <-> ((z ε x) & (z ε y)))
 AndInt 20 37
\end{verbatim}	

with log

\begin{verbatim}
21. Hyp("Elem(z, intersection(x,y))")
22. DefEqInt(1)
23. EqualitySub(21,22,[0])
24. ClassElim(23)
25. AndElimR(24)
26. ImpInt(25,21)
27. Hyp("(Elem(z,x) & Elem(z,y))")
28. AndElimL(27)
29. ExistsInt(28,"x","x",[0])
30. DefSub(29,"Set",["z"],[0])
31. AndInt(30,27)
32. ClassInt(31,"z")
33. Symmetry(22)
34. EqualitySub(32,33,[0])
35. ImpInt(34,27)
36. AndInt(26,35)
37. EquivConst(36)
38. AndInt(20,37)
\end{verbatim}

This theorem comes in the main directory of the PyLog repository and can be loaded via Load("Th4"). We also show the proof of  the first half of the conjunction of theorem 5 of \cite{kel} which illustrates how we use previous theorems and apply an axiom:

\begin{verbatim}
0. z ε (x ∪ x) Hyp 
1. ((z ε (x ∪ y)) <-> ((z ε x) v (z ε y))) & ((z ε (x ∩ y)) <-> ((z ε x) & (z ε y))) TheoremInt 
2. (z ε (x ∪ y)) <-> ((z ε x) v (z ε y)) AndElimL 1
3. ((z ε (x ∪ y)) -> ((z ε x) v (z ε y))) & (((z ε x) v (z ε y)) -> (z ε (x ∪ y))) EquivExp 
4. (z ε (x ∪ y)) -> ((z ε x) v (z ε y)) AndElimL 3
5. ∀y.((z ε (x ∪ y)) -> ((z ε x) v (z ε y))) ForallInt 4
6. (z ε (x ∪ x)) -> ((z ε x) v (z ε x)) ForallElim 5
7. (z ε x) v (z ε x) ImpElim 0 6
8. z ε x Hyp 
9. z ε x Hyp 
10. z ε x OrElim 7 8 8 9 9
11. (z ε (x ∪ x)) -> (z ε x) ImpInt 10 Qed
12. z ε x Hyp 
13. (z ε x) v (z ε x) OrIntL 12
14. ((z ε x) v (z ε y)) -> (z ε (x ∪ y)) AndElimR 3
15. ∀y.(((z ε x) v (z ε y)) -> (z ε (x ∪ y))) ForallInt 14
16. ((z ε x) v (z ε x)) -> (z ε (x ∪ x)) ForallElim 15
17. z ε (x ∪ x) ImpElim 13 16
18. (z ε x) -> (z ε (x ∪ x)) ImpInt 17 Qed
19. ((z ε (x ∪ x)) -> (z ε x)) & ((z ε x) -> (z ε (x ∪ x))) AndInt 11 18
20. (z ε (x ∪ x)) <-> (z ε x) EquivConst 
21. ∀z.((z ε (x ∪ x)) <-> (z ε x)) ForallInt 20 Qed
22. ∀x.∀y.((x = y) <-> ∀z.((z ε x) <-> (z ε y))) AxInt 
23. ∀y.(((x ∪ x) = y) <-> ∀z.((z ε (x ∪ x)) <-> (z ε y))) ForallElim 22
24. ((x ∪ x) = x) <-> ∀z.((z ε (x ∪ x)) <-> (z ε x)) ForallElim 23
25. (((x ∪ x) = x) -> ∀z.((z ε (x ∪ x)) <-> (z ε x))) & (∀z.((z ε (x ∪ x)) <-> (z ε x)) -> ((x ∪ x) = x)) EquivExp 
26. ∀z.((z ε (x ∪ x)) <-> (z ε x)) -> ((x ∪ x) = x) AndElimR 25
27. (x ∪ x) = x ImpElim 21 26 Qed
\end{verbatim}

the log is

\begin{verbatim}
0. Hyp("Elem(z,union(x,x))")
1. TheoremInt(1)
2. AndElimL(1)
3. EquivExp(2)
4. AndElimL(3)
5. ForallInt(4,"y","y")
6. ForallElim(5,"x")
7. ImpElim(0,6)
8. Hyp("Elem(z,x)")
9. Hyp("Elem(z,x)")
10. OrElim(7,8,8,9,9)
11. ImpInt(10,0)
12. Hyp("Elem(z,x)")
13. OrIntL(12,"Elem(z,x)")
14. AndElimR(3)
15. ForallInt(14,"y","y")
16. ForallElim(15,"x")
17. ImpElim(13,16)
18. ImpInt(17,12)
19. AndInt(11,18)
20. EquivConst(19)
21. ForallInt(20,"z","z")
22. AxInt(0)
23. ForallElim(22,"union(x,x)")
24. ForallElim(23,"x")
25. EquivExp(24)
26. AndElimR(25)
27. ImpElim(21,26)
\end{verbatim}

To prove Th6 and Th7 we could make use of PolySub and previously proven propositional validities.

For instance if a logical validity was previously proven

\begin{verbatim}
0. A v B Hyp 
1. A Hyp 
2. B v A OrIntL 1
3. B Hyp 
4. B v A OrIntR 3
5. B v A OrElim 0 1 2 3 4
6. (A v B) -> (B v A) ImpInt 5 Qed

0. Hyp("(A v B)")
1. Hyp("A")
2. OrIntL(1,"B")
3. Hyp("B")
4. OrIntR(3,"A")
5. OrElim(0,1,2,3,4)
6. ImpInt(5,0)
\end{verbatim}

Suppose we saved this theorem as "Log1".
Then we can add this theorem to our environment
and instantiate it via PolySub to, for instance:

\begin{verbatim}
>>> Load("Kelley-Morse")
True
>>> 
>>> 
>>> ViewTheorem("Log1")
'Log1 : (A v B) -> (B v A)'
>>> LoadTheorem("Log1")
True
>>> ShowTheorems()
0. A v ¬A 
1. (A v B) -> (B v A) 
>>> TheoremInt(1)
True
>>> ShowProof()
0. (A v B) -> (B v A) TheoremInt 
>>> PolySub(0,"A","Elem(z,x)")
(...)
>>> PolySub(1,"B","Elem(z,y)")
0. (A v B) -> (B v A) TheoremInt 
1. ((z ε x) v B) -> (B v (z ε x)) PolySub 0
2. ((z ε x) v (z ε y)) -> ((z ε y) v (z ε x)) PolySub 1
>>> ShowLog()
0. TheoremInt(1)
1. PolySub(0,"A","Elem(z,x)")
2. PolySub(1,"B","Elem(z,y)")
\end{verbatim}
 It is important that when adding theorems the theorems were saved with the same environment. An exception occurs for logical validities.

If a predicate is part of the language but has not been defined then we can consider it a second-order variable in the classical sense (i.e. having a fixed arity) and
use the command PredSub. This can  be useful for implementing axiom schemes such as induction.  Also it is useful for reusing first-order validities (which is not possible with the polymorphic variables).

\begin{verbatim}
	
>>> Load("Kelley-Morse")
True
>>> AddPredicate("Foo",1,True)
True
>>> Hyp("Foo(x)")
0. Foo(x)  Hyp 
True
>>> PredSub(0,"Foo",["x"], "neg Set(x)", [0])
0. Foo(x)  Hyp 
1. ¬Set(x)  PredSub 0
True
>>> PredSub(1,"Set",["x"], "neg Set(x)", [0])
Predicate is defined.	
	
\end{verbatim}

To test a theory, a collection of theorems, one can run the CheckTheory function,
For instance, for the theorems up to theorem 55 in Kelley-Morse  set theory:

\begin{verbatim}
>>> CheckTheory(["Th4","Th5","Th6","Th7","Th8",
"Th11","Th12","Th14","Th16","Th17","Th19","Th20",
"Th21","Th24","Th26","Th27","Th28","Th29","Th30",
"Th31","Th32","Th33", "Th34","Th35", "Th37","Th38",
"Th39","Th41", "Th42","Th43","Th44","Th46","Th47","Th49","Th50","Th53","Th54","Th55", "Th58", "Th59", "Th61", "Th62", "Th64", "Th67", "Th69", "Th70", "Th71", "Th73", "Th74","Th75","Th77"])
\end{verbatim}

This will generate each proof again from the log. If the final line has Qed then the check is succesful.

\section*{Proof of Law of Excluded Middle}

We present here a proof in PyLog (using AbsC )of the useful classical validity $A v \neg A$. To do this we use
three validities:

\begin{verbatim}
>>> ShowTheorems()
0. D <-> ¬¬D 
1. (A -> B) -> (¬B -> ¬A) 
2. ¬(A v B) <-> (¬A & ¬B) 
\end{verbatim}

The proofs of these theorems is an easy exercise. Only theorem 0 uses AbsC. Here is the proof saved under the name "ExcludedMiddle".

\begin{verbatim}
>>> ShowProof()
0. ¬(A v B) <-> (¬A & ¬B) TheoremInt 
1. ¬(A v ¬A) <-> (¬A & ¬¬A) PolySub 0
2. (¬(A v ¬A) -> (¬A & ¬¬A)) & ((¬A & ¬¬A) -> ¬(A v ¬A)) EquivExp 
3. D <-> ¬¬D TheoremInt 
4. (D -> ¬¬D) & (¬¬D -> D) EquivExp 
5. ¬A & ¬¬A Hyp 
6. ¬A AndElimL 5
7. ¬¬A AndElimR 5
8. ¬¬D -> D AndElimR 4
9. ¬¬A -> A PolySub 8
10. A ImpElim 7 9
11. ¬A & A AndInt 6 10
12. (¬A & ¬¬A) -> (¬A & A) ImpInt 11 Qed
13. ¬(A v ¬A) -> (¬A & ¬¬A) AndElimL 2
14. ¬(A v ¬A) Hyp 
15. ¬A & ¬¬A ImpElim 14 13
16. ¬A & A ImpElim 15 12
17. ¬(A v ¬A) -> (¬A & A) ImpInt 16 Qed
18. (A -> B) -> (¬B -> ¬A) TheoremInt 
19. (¬(A v ¬A) -> B) -> (¬B -> ¬¬(A v ¬A)) PolySub 18
20. (¬(A v ¬A) -> (¬A & A)) -> (¬(¬A & A) -> ¬¬(A v ¬A)) PolySub 19
21. ¬(¬A & A) -> ¬¬(A v ¬A) ImpElim 17 20
22. ¬¬(A v ¬A) -> (A v ¬A) PolySub 8
23. ¬(¬A & A) Hyp 
24. ¬¬(A v ¬A) ImpElim 23 21
25. A v ¬A ImpElim 24 22
26. ¬(¬A & A) -> (A v ¬A) ImpInt 25 Qed
27. ¬A & A Hyp 
28. ¬A AndElimL 27
29. A AndElimR 27
30. _|_ ImpElim 29 28
31. ¬(¬A & A) ImpInt 30 Qed
32. A v ¬A ImpElim 31 26 Qed
>>> ShowLog()
0. TheoremInt(2)
1. PolySub(0,"B","neg A")
2. EquivExp(1)
3. TheoremInt(0)
4. EquivExp(3)
5. Hyp("(neg A & neg neg A)")
6. AndElimL(5)
7. AndElimR(5)
8. AndElimR(4)
9. PolySub(8,"D","A")
10. ImpElim(7,9)
11. AndInt(6,10)
12. ImpInt(11,5)
13. AndElimL(2)
14. Hyp("neg (A v neg A)")
15. ImpElim(14,13)
16. ImpElim(15,12)
17. ImpInt(16,14)
18. TheoremInt(1)
19. PolySub(18,"A","neg (A v neg A)")
20. PolySub(19,"B","(neg A & A)")
21. ImpElim(17,20)
22. PolySub(8,"D","(A v neg A)")
23. Hyp("neg (neg A & A)")
24. ImpElim(23,21)
25. ImpElim(24,22)
26. ImpInt(25,23)
27. Hyp("(neg A & A)")
28. AndElimL(27)
29. AndElimR(27)
30. ImpElim(29,28)
31. ImpInt(30,27)
32. ImpElim(31,26)
\end{verbatim}

\section*{The Gentzen Automatic Theorem Prover}

Included within PyLog is a simple algorithm for finding proofs of propositiona validities in the intuitionistic propositional calculus. Smullyan \cite{smu} has remarked that tableaux for classical propositional calculus correspond to proofs in the sequent calculus "turned upside-down".
This is a paradigm of automatic proof generation.
The \emph{cut-free} sequent calculi for the intuitionistic propositional calculus (IPC), more specifically, the Gentzen system G3i
restricted to IPC (wherein the usual structural rules are "absorbed") is a striking example of this correspondence. 
Indeed, it is well known that G3i restricted to propositional formulas can be "inverted" to furnish a 
decision procedure for IPC \cite{tro}[4.2.6].
For propositional G3i we have the rules:

\begin{center}
	\AxiomC{}
	\RightLabel{Ax}
	\UnaryInfC{$P,\Gamma\Rightarrow P$}
	\DisplayProof
	\hskip 6.5em
	\AxiomC{}
	\RightLabel{L$\bot$}
	\UnaryInfC{$\bot,\Gamma\Rightarrow A $}
	\DisplayProof
\end{center}

\begin{center}
	\AxiomC{$A,B,\Gamma\Rightarrow C$}
	\RightLabel{L$\wedge$}
	\UnaryInfC{$A\wedge B,\Gamma\Rightarrow C$}
	\DisplayProof
	\hskip 6.5em
	\AxiomC{$\Gamma\Rightarrow A$}
	\AxiomC{$\Gamma\Rightarrow B$}
	\RightLabel{R$\wedge$}
	\BinaryInfC{$\Gamma\Rightarrow A\wedge B$}
	\DisplayProof
\end{center}

\begin{center}
	\AxiomC{$A,\Gamma\Rightarrow C$}
	\AxiomC{$B,\Gamma\Rightarrow C$}
	\RightLabel{L$\vee$}
	\BinaryInfC{$A\vee B,\Gamma\Rightarrow C$}
	\DisplayProof
	\hskip 6.5em
	\AxiomC{$\Gamma\Rightarrow A_i$}
	\RightLabel{R$\vee, i=0,1$}
	\UnaryInfC{$\Gamma\Rightarrow A_0\vee A_1$}
	\DisplayProof
\end{center}

\begin{center}
	\AxiomC{$A\rightarrow B,\Gamma\Rightarrow A$}
	\AxiomC{$B,\Gamma\Rightarrow C$}
	\RightLabel{L$\rightarrow$}
	\BinaryInfC{$A\rightarrow B,\Gamma\Rightarrow C$}
	\DisplayProof
	\hskip 6.5em
	\AxiomC{$A,\Gamma\Rightarrow B$}
	\RightLabel{R$\rightarrow$}
	\UnaryInfC{$\Gamma\Rightarrow A\rightarrow B$}
	\DisplayProof
\end{center}

The rules either have no premises (axioms) or else one or two premises.
To turn a proof upside-down we obviously have to somehow "invert" the rules.  For
a special class of such systems (those that are \emph{length-preserving} and \emph{analytic}, which corresponds to
the sub-formula property) we can perform an "inversion" to obtain what we call a \emph{reductive system}
which is easier to work with for the generation of proofs. A proof is seen as a winning \emph{linear sequence
	of moves} which when acting upon a stateful system yields the empty set.

Our algorithm adapts the proof of \cite{tro}[4.2.6]. The basic object is a \emph{sequence list object} $S$ which includes a list of sequents $(s_1,...,s_n)$  and a state $(m,h)$ which includes the history $h$ of the rules that have been applied to it and the accumulated set  $m$ of all the new sequents that have been generated by these rules (this is to avoid cycles). A rule $r$  is an application which takes a sequence list $S$ and yields a new sequence list $S'$. The transformation of the list of sequence is effected by choosing one specific sequence $s_k$, eliminating it and possibly generating new sequents  which are incorporate into the previous list with $s_k$ removed. The rule is added to $h$ and the new sequents (if any) are added to $m$.
To prove a formula $f$ we start with a sequence object $(\Rightarrow f), (\emptyset,()))$, that is, having a single sequent. We want to find a sequence of rules $r_1,...,r_n$ such that when applied in order yield a sequent object of the form $((), (m,h))$.
The best way to understand our inversion of propositional G3i is to examine the source code which is found in the file gentzen.py in the PyLog main folder. Consider the Sequent and SequentList classes:

\begin{verbatim}
class Sequent:
  def __init__(self, head,body):
   self.head = head
   self.body = body
   
(...)   

class SequentList:
 def __init__(self,formstring):
  form = astop.NegationExpand(parser.Formula(tokenizer.Tokenize(formstring)))
  seq = Sequent(form,[])
  self.sequentlist = [seq]
  self.memory = []
  self.proof = []

(...)
\end{verbatim}

Here is an example of a rule:

\begin{verbatim}
def rand(seq,n):
 if  n < len(seq.sequentlist)  and seq.sequentlist[n].head.name=="constructor":
  if seq.sequentlist[n].head.operator.name=="&":
   aux1 = seq.sequentlist[:n]
   aux2 = seq.sequentlist[n+1:]
   aux3 = seq.sequentlist[n].body
   left = seq.sequentlist[n].head.left
   right = seq.sequentlist[n].head.right
   new1 = Sequent(left,aux3)
   new2 = Sequent(right, aux3)
   if not SeqInc(new1, seq.memory) and not SeqInc(new2,seq.memory):
    seq.sequentlist = aux1 + [new1,new2] + aux2
    seq.memory.append(new1)
    seq.memory.append(new2)
    seq.proof.append("rand(" + str(n) +")")
   return seq
 return False
\end{verbatim}

There are nine rules rand, ror1,ror2, rimp, land, lor, limp, labs, ax.

This genzen.py program runs on the same syntactic engine as PyLog and formulas are entered the same way. Note that it can only find proofs of \emph{intuistionistic} propositional validities.
The  algorithm  used to find a proof is a simple brute-force search which avoids cycles. It can take some time.

\begin{verbatim}
$ python3.9 -i gentzen.py
>>> Auto("( neg (A v B) -> (neg A & neg B))")
>>> Prove()
1. rimp(0)
2. rand(0)
3. rimp(0)
4. imp(0,1)
5. ror1(0)
6. ax(0,0)   A, ¬(A v B) => A
7. ax(0,1)   A, _|_ => _|_
8. rimp(0)
9. limp(0,1)
10. ror2(0)
11. ax(0,0)   B, ¬(A v B) => B
12. ax(0,1)   B, _|_ => _|_

0.  B, ¬(A v B) => B   Ax0  
1.  B, _|_ => _|_   Ax1  
2.  B, ¬(A v B) => A v B   Ror2  0 
3.  B, ¬(A v B) => _|_   Limp1  2 1 
4.  A, ¬(A v B) => A   Ax0  
5.  A, _|_ => _|_   Ax1  
6.  A, ¬(A v B) => A v B   Ror1  4 
7.  A, ¬(A v B) => _|_   Limp1  6 5 
8.  ¬(A v B) => ¬B   Rimp  3 
9.  ¬(A v B) => ¬A   Rimp  7 
10.  ¬(A v B) => ¬A & ¬B   Rand  9 8 
11.  => ¬(A v B) -> (¬A & ¬B)   Rimp  10 
\end{verbatim}

As seen above there  is also an algorithm which reconstructs a linear version of the full sequent proof tree from the proof. Looking at this
linear sequent proof it is very easy to obtain the corresponding linear natural deduction proof. This will be implemented in future versions of Pylog.
After the  proof sequence  is obtained can also be tested manually by applying the rules in order to the initial
sequent list object (which is the object State).

\begin{verbatim}
>>> State.display()
0.  => ¬(A v B) -> (¬A & ¬B)
>>> rimp(State,0)
<__main__.SequentList object at 0x109a69070>
>>> State.display()
0.  ¬(A v B) => ¬A & ¬B
>>> rand(State,0)
<__main__.SequentList object at 0x109a69070>
>>> State.display()
0.  ¬(A v B) => ¬A
1.  ¬(A v B) => ¬B
>>> rimp(State,0)
<__main__.SequentList object at 0x109a69070>
>>> State.display()
0.  A, ¬(A v B) => _|_
1.  ¬(A v B) => ¬B
>>> limp(State,0,1)
<__main__.SequentList object at 0x109a69070>
>>> State.display()
0.  A, ¬(A v B) => A v B
1.  A, _|_ => _|_
2.  ¬(A v B) => ¬B
\end{verbatim}

and so forth. We will in the future implement a further algorithm that transforms the linear sequent proof tree into a Pylog linear natural deduction proof.

\section*{Formalising Category Theory}

\begin{verbatim}
>>> ShowAxioms()
0. ([f: A → B|C] & [g: B → D|C]) -> [(g∘f): A → D|C] 
1. A:C -> [id(A): A → A|C] 
2. (A:C & [f: A → B|C]) -> (((id(A)∘f) = f) & ((f∘id(A)) = f)) 
3. ([f: A → B|C] & ([g: B → D|C] & [h: D → E|C])) -> ((h∘(g∘f)) = ((h∘g)∘f)) 
4. A:C -> Cat(C) 
5. [f: A → B|C] -> (A:C & B:C) 
6. (Cat(A) & Cat(B)) <-> Cat(func(A,B)) 
7. F:func(A,B) -> FA:B 
8. (F:func(A,B) & [f: a → b|A]) -> [Ff: Fa → Fb|B] 
9. (F:func(A,B) & a:A) -> (Fid(A) = id(FA)) 
10. (F:func(A,B) & ([f: a → b|C] & [g: b → c|C])) -> (F(g∘f) = (Fg∘Ff)) 
11. ([e: F → G|func(A,B)] & a:A) -> [nat(e,a): Fa → Ga|B] 
12. ([e: F → G|func(A,B)] & [f: a → b|A]) -> ((Gf∘nat(e,a)) = (nat(e,b)∘Ff)) 


>>> ShowDefinitions()
Terminal(A,C) <-> (A:C & ∀B.(B:C -> ∃¹f.[f: B → A|C]))
Isomorphism(f,A,B,C) <-> ([f: A → B|C] & ∃g.([g: B → A|C] & (((g∘f) = id(A)) & ((f∘g) = id(B)))))
Iso(A,B,C) <-> ∃f.Isomorphism(f,A,B,C)

>>> ShowProof()
0. Terminal(T,C)  Hyp 
1. Terminal(S,C)  Hyp 
2. T:C & ∀B.(B:C -> ∃¹f.[f: B → T|C])  DefExp 0
3. S:C & ∀B.(B:C -> ∃¹f.[f: B → S|C])  DefExp 1
4. T:C  AndElimL 2
5. S:C  AndElimL 3
6. ∀B.(B:C -> ∃¹f.[f: B → T|C])  AndElimR 2
7. ∀B.(B:C -> ∃¹f.[f: B → S|C])  AndElimR 3
8. T:C -> ∃¹f.[f: T → T|C]  ForallElim 6
9. S:C -> ∃¹f.[f: S → S|C]  ForallElim 7
10. A:C -> [id(A): A → A|C]  AxInt 
11. ∀A.(A:C -> [id(A): A → A|C])  ForallInt 10
12. T:C -> [id(T): T → T|C]  ForallElim 11
13. ∀A.(A:C -> [id(A): A → A|C])  ForallInt 10
14. S:C -> [id(S): S → S|C]  ForallElim 13
15. [id(T): T → T|C]  ImpElim 4 12
16. [id(S): S → S|C]  ImpElim 5 14
17. ∃¹f.[f: T → T|C]  ImpElim 4 8
18. ∃¹f.[f: S → S|C]  ImpElim 5 9
19. ∃f.([f: T → T|C] & ∀l.([l: T → T|C] -> (l = f)))  UniqueElim 17
20. ∃f.([f: S → S|C] & ∀m.([m: S → S|C] -> (m = f)))  UniqueElim 18
21. [j: T → T|C] & ∀l.([l: T → T|C] -> (l = j))  Hyp 
22. [k: S → S|C] & ∀m.([m: S → S|C] -> (m = k))  Hyp 
23. ∀l.([l: T → T|C] -> (l = j))  AndElimR 21
24. ∀m.([m: S → S|C] -> (m = k))  AndElimR 22
25. [j: T → T|C]  AndElimL 21
26. [k: S → S|C]  AndElimL 22
27. S:C -> ∃¹f.[f: S → T|C]  ForallElim 6
28. T:C -> ∃¹f.[f: T → S|C]  ForallElim 7
29. ∃¹f.[f: S → T|C]  ImpElim 5 27
30. ∃¹f.[f: T → S|C]  ImpElim 4 28
31. ∃f.([f: S → T|C] & ∀r.([r: S → T|C] -> (r = f)))  UniqueElim 29
32. ∃f.([f: T → S|C] & ∀s.([s: T → S|C] -> (s = f)))  UniqueElim 30
33. [u: S → T|C] & ∀r.([r: S → T|C] -> (r = u))  Hyp 
34. [v: T → S|C] & ∀s.([s: T → S|C] -> (s = v))  Hyp 
35. [u: S → T|C]  AndElimL 33
36. [v: T → S|C]  AndElimL 34
37. ([f: A → B|C] & [g: B → D|C]) -> [(g∘f): A → D|C]  AxInt 
38. ∀A.(([f: A → B|C] & [g: B → D|C]) -> [(g∘f): A → D|C])  ForallInt 37
39. ([f: T → B|C] & [g: B → D|C]) -> [(g∘f): T → D|C]  ForallElim 38
40. ∀B.(([f: T → B|C] & [g: B → D|C]) -> [(g∘f): T → D|C])  ForallInt 39
41. ([f: T → S|C] & [g: S → D|C]) -> [(g∘f): T → D|C]  ForallElim 40
42. ∀f.(([f: T → S|C] & [g: S → D|C]) -> [(g∘f): T → D|C])  ForallInt 41
43. ([v: T → S|C] & [g: S → D|C]) -> [(g∘v): T → D|C]  ForallElim 42
44. ∀g.(([v: T → S|C] & [g: S → D|C]) -> [(g∘v): T → D|C])  ForallInt 43
45. ([v: T → S|C] & [u: S → D|C]) -> [(u∘v): T → D|C]  ForallElim 44
46. ∀D.(([v: T → S|C] & [u: S → D|C]) -> [(u∘v): T → D|C])  ForallInt 45
47. ([v: T → S|C] & [u: S → T|C]) -> [(u∘v): T → T|C]  ForallElim 46
48. [v: T → S|C] & [u: S → T|C]  AndInt 36 35
49. [(u∘v): T → T|C]  ImpElim 48 47
50. ∀A.(([f: A → B|C] & [g: B → D|C]) -> [(g∘f): A → D|C])  ForallInt 37
51. ([f: S → B|C] & [g: B → D|C]) -> [(g∘f): S → D|C]  ForallElim 50
52. ∀B.(([f: S → B|C] & [g: B → D|C]) -> [(g∘f): S → D|C])  ForallInt 51
53. ([f: S → T|C] & [g: T → D|C]) -> [(g∘f): S → D|C]  ForallElim 52
54. ∀D.(([f: S → T|C] & [g: T → D|C]) -> [(g∘f): S → D|C])  ForallInt 53
55. ([f: S → T|C] & [g: T → S|C]) -> [(g∘f): S → S|C]  ForallElim 54
56. ∀f.(([f: S → T|C] & [g: T → S|C]) -> [(g∘f): S → S|C])  ForallInt 55
57. ([u: S → T|C] & [g: T → S|C]) -> [(g∘u): S → S|C]  ForallElim 56
58. ∀g.(([u: S → T|C] & [g: T → S|C]) -> [(g∘u): S → S|C])  ForallInt 57
59. ([u: S → T|C] & [v: T → S|C]) -> [(v∘u): S → S|C]  ForallElim 58
60. [u: S → T|C] & [v: T → S|C]  AndInt 35 36
61. [(v∘u): S → S|C]  ImpElim 60 59
62. [id(T): T → T|C] -> (id(T) = j)  ForallElim 23
63. id(T) = j  ImpElim 15 62
64. [(u∘v): T → T|C] -> ((u∘v) = j)  ForallElim 23
65. (u∘v) = j  ImpElim 49 64
66. [id(S): S → S|C] -> (id(S) = k)  ForallElim 24
67. id(S) = k  ImpElim 16 66
68. [(v∘u): S → S|C] -> ((v∘u) = k)  ForallElim 24
69. (v∘u) = k  ImpElim 61 68
70. j = id(T)  Symmetry 63
71. k = id(S)  Symmetry 67
72. (u∘v) = id(T)  EqualitySub 65 70
73. (v∘u) = id(S)  EqualitySub 69 71
74. ((u∘v) = id(T)) & ((v∘u) = id(S))  AndInt 72 73
75. [u: S → T|C] & (((u∘v) = id(T)) & ((v∘u) = id(S)))  AndInt 35 74
76. ∃h.([h: S → T|C] & (((h∘v) = id(T)) & ((v∘h) = id(S))))  ExistsInt 75
77. [v: T → S|C] & ∃h.([h: S → T|C] & (((h∘v) = id(T)) & ((v∘h) = id(S))))  AndInt 36 76
78. Isomorphism(v,T,S,C)  DefSub 77
79. ∃f.Isomorphism(f,T,S,C)  ExistsInt 78
80. Iso(T,S,C)  DefSub 79
81. Iso(T,S,C)  ExistsElim 19 21 80
82. Iso(T,S,C)  ExistsElim 20 22 81
83. Iso(T,S,C)  ExistsElim 31 33 82
84. Iso(T,S,C)  ExistsElim 32 34 83
85. Terminal(S,C) -> Iso(T,S,C)  ImpInt 84
86. Terminal(T,C) -> (Terminal(S,C) -> Iso(T,S,C))  ImpInt 85 Qed
	
	
	
\end{verbatim}

\section*{Principle of Transcendental Analogy}

It is cumbersome to work with formal proofs involving quotients of algebraic structures.
We consider a class $U$ of all the different models (or representatives of isomorphism classes) of a given structure. And we fix an axiomatic system for operations $f_1,...,f_n$ which act upon the whole $\bigcup U$ and by restrictions on all the elements $M$ of  $U$ where it corresponds to the various operations of these structures.  It is like Weil's foundations for Algebraic Geometry in which all fields of rational functions are seen as subfields of one large field.
Quotients (for instance for groups or rings) are defined not in terms of structures but of epimorphisms $f:A\rightarrow B$  and a given substructure (normal subgroup, ideal) of the domain $A$.
\section*{Second-Order Arithmetic}

\begin{verbatim}
The project is to formalise Stephen G. Simpson's book Subsystems of 
Second-Order Arithmetic. The axioms of Z2 are interpreted as follows:

0. Nat(x) v Set(x) 
1. Nat(x) -> ¬Set(x) 
2. Nat(0) & Nat(1) 
3. Nat(n) -> ¬((n + 1) = 0) 
4. (Nat(n) & Nat(m)) -> (((m + 1) = (n + 1)) -> (m = n)) 
5. Nat(m) -> ((m + 0) = m) 
6. (Nat(n) & Nat(m)) -> ((m + (n + 1)) = ((m + n) + 1)) 
7. Nat(m) -> ((m * 0) = 0) 
8. (Nat(n) & Nat(m)) -> ((m * (n + 1)) = ((m * n) + m)) 
9. Nat(m) -> ¬(m < 0) 
10. (Nat(n) & Nat(m)) -> ((m < (n + 1)) <-> ((m < n) v (m = n))) 
11. (x ε X) -> (Nat(x) & Set(X)) 
12. ((0 ε X) & ∀n.((n ε X) -> ((n + 1) ε X))) -> ∀n.(n ε X) 

n ε { x: form(x)}
-----------------
Nat(n) & form(n)


Nat(n) & form(n)
-----------------
n ε { x: form(x)}

The comprehension scheme is converted into the two rules above which are similar
to ClassElim and ClassInt used to formalise Kelley-Morse
set theory. Example of a simple proof:

0. 1 = 0  Hyp 
1. Nat(0) & Nat(1)  AxInt 
2. 0 = 0  Identity 
3. Nat(0)  AndElimL 1
4. Nat(n) -> ¬((n + 1) = 0)  AxInt 
5. ∀n.(Nat(n) -> ¬((n + 1) = 0))  ForallInt 4
6. Nat(0) -> ¬((0 + 1) = 0)  ForallElim 5
7. ¬((0 + 1) = 0)  ImpElim 3 6
8. ¬((0 + 0) = 0)  EqualitySub 7 0
9. Nat(m) -> ((m + 0) = m)  AxInt 
10. ∀m.(Nat(m) -> ((m + 0) = m))  ForallInt 9
11. Nat(0) -> ((0 + 0) = 0)  ForallElim 10
12. (0 + 0) = 0  ImpElim 3 11
13. ¬(0 = 0)  EqualitySub 8 12
14. _|_  ImpElim 2 13
15. ¬(1 = 0)  ImpInt 14 Qed

\end{verbatim}

\section*{Euclidean Geometry}

For a formalisation of the first theorem of Euclid in PyLog see \cite{eucl}.


\begin{thebibliography}{1}
		\bibitem{eucl} Protin, C. L. (2022). \emph{Euclidean Logic}. \url{http://dx.doi.org/10.13140/RG.2.2.15673.65120}.
	\bibitem{eucl} Protin, C. L. (2022). \emph{Kelley-Morse Set Theory in PyLog}. \url{http://dx.doi.org/10.13140/RG.2.2.13288.14081}.
	\bibitem{dag} D. Prawitz, \emph{Natural Deduction}, Dover.
	\bibitem{kel} J. Kelley, \emph{General Topology}.
	\bibitem{tro} A. Troelstra, \emph{Constructivism in Mathematics} vol. I.	
	\bibitem{smu}
	Raymond M. Smullyan,
	\textit{First-Order Logic}, Dover 1995.
	
	\bibitem{gir} 
	J.-Y. Girard,
	\textit{Proofs and Types}, Cambridge University Press, 1989.
	
	
	\bibitem{schwich} 
	A.S. Troelstra, H. Schwichtenberg,
	\textit{Basic Proof Theory}, Cambridge Tracts on Theoretical Computer Science, 2nd Ed., 2000.
	
	
\end{thebibliography}
 \end{document}